\documentclass[11pt]{amsart}

\usepackage{amsmath,amsthm,amssymb,amssymb, paralist, xspace, graphicx, url, amscd, euscript, mathrsfs,stmaryrd,epic,eepic,color}

\usepackage[all]{xy}
\usepackage{amsmath,amscd}
\SelectTips{cm}{}

\usepackage[colorlinks=true]{hyperref}

\usepackage{pstricks}

\numberwithin{equation}{subsection}

1

\numberwithin{subsection}{section}

\allowdisplaybreaks[1]

\newtheorem*{namedtheorem}{\theoremname}
\newcommand{\theoremname}{testing}

\theoremstyle{plain}

\newtheorem{thm}{Theorem}[section]
\newtheorem{proposition}[thm]{Proposition}
\newtheorem{proposition-definition}[thm]{Proposition-Definition}
\newtheorem{lemma-definition}[thm]{Lemma-Definition}
\newtheorem{corollary}[thm]{Corollary}
\newtheorem{lemma}[thm]{Lemma}
\newtheorem{conj}[thm]{Conjecture}

\theoremstyle{definition}
\newtheorem{definition}[thm]{Definition}
\newtheorem{notation}[thm]{Notation}

\newtheorem{remark}[thm]{Remark}

\theoremstyle{remark}

\numberwithin{thm}{section}

\newcommand\cA{\mathcal{A}}

\newcommand\cM{\mathcal{M}}

\newcommand\N{\mathcal{N}}

\def\O{\mathcal{O}}

\def\P{\mathbb{P}}
\def\A{\mathbb{A}}
\def\N{\mathbb{N}}

\def\C{\mathbb{C}}

\def\im{Im}
\def\tilde{\widetilde}

\newcommand{\R}{}

\newcommand\arr{\ifinner\to\else\longrightarrow\fi}

\def\displaytimes_#1{\mathrel{\mathop{\times}\limits_{#1}}}

\def\displayotimes_#1{\mathrel{\mathop{\bigotimes}\limits_{#1}}}

\newcommand\CH{\operatorname{CH}}

\newcommand\Spec{\operatorname{Spec}}

\newcommand\sym{\operatorname{Sym}}

\newdir{ >}{{}*!/-5pt/@{>}}

\newcommand\doublelong[2]{\mathbin{\xymatrix{{}\ar@<3pt>[r]^{#1}
\ar@<-3pt>[r]_{#2}&}}}

\newlength{\ignora}

\newcommand{\Z}{\mathbb{Z}}

\numberwithin{equation}{subsection}

\newcommand{\Alb}{\operatorname{Alb}}

\begin{document}


\title{$\A^1$-equivalence of zero cycles on surfaces}


\author{Yi Zhu}

\address[Zhu]{Pure Mathematics\\Univeristy of Waterloo\\Waterloo, ON N2L3G1\\ Canada}
\email{yi.zhu@uwaterloo.ca}



\begin{abstract}
In this paper, we study $\A^1$-equivalence classes of zero cycles on open algebraic surfaces. We prove the logarithmic version of Mumford's theorem on zero cycles. We also prove that the log Bloch conjecture holds for  surfaces with log Kodaira dimension $-\infty$.

\end{abstract}
\maketitle


\section{Introduction}\label{sec:intro}
Let $X$ be a smooth projective complex surface. Understanding the structure of the Chow group of zero cycles of degree zero $\CH_0(X)^0$ is important but difficult. Mumford first studied this group and proved the following theorem.
\begin{thm}\cite{Mumford68}
If $h^0(X,\Omega^2_X)> 0$, the group $\CH_0(X)^0$ is infinite-dimensional.
\end{thm}

In the other direction, we have Bloch's conjecture as below. 
\begin{conj}\cite{Bloch80}
If $h^0(X,\Omega^2_X)= 0$, then the Albanese morphism induces an isomorphism $$\CH_0(X)^0\cong\Alb(X).$$
\end{conj}

Bloch's conjecture has been proved for smooth projective surfaces with Kodaira dimension less than two \cite{Bloch76}. For surfaces of general type, many cases have been proved, but it is still widely open in general \cite[Chapter 11]{Voisin}.

For not necessarily proper varieties, Spie{\ss} and Szamuely \cite{Szamuely03} observe that the right replacement for Chow group of zero cycles is Suslin's $0$-th algebraic singular homology $h_0(U)^0$ and furthermore they prove the log Roitmann's theorem for smooth quasiprojective varieties in all dimensions.

\begin{definition}
	Let $U$ be a smooth quasiprojective variety. Two zero cycles $A_1$, $A_2$ of degree $n$ are \emph{$\A^1$-equivalent} if there exists a $0$-cycle $B$ of degree $m$ such that
	\begin{itemize}
		\item $A_1+B$ and $A_2+B$ are effective;
		\item there exists a morphism $z:\A^1\to \sym^{n+m}U$ such that $z(0)=A_1+B$, $z(1)=A_2+B$.
	\end{itemize} 
\end{definition}

\begin{definition}
	Suslin's zeroth homology $h_0(U)^0$ is the group of all zero-cycles on $U$ of degree $0$ modulo $\A^1$-equivalences. 
\end{definition}

When $U$ is a curve, $\A^1$-equivalence is indeed the equivalence relation of divisors defined by the modulus $D$ as in \cite[V.2]{Serre-AG}.

\begin{thm}\cite[Theorem 1.1]{Szamuely03}
Given a smooth quasiprojective variety $U$, the Albanese morphism 
$$alb:h_0(U)^0\to \Alb(U)$$
induces an isomorphism on the torsion subgroups. 
\end{thm}

It is natural to consider Mumford's theorem for smooth quasiprojective surfaces. Consider the map$$\sigma_d:Sym^d(U)\times Sym^d(U)\to h_0(U)^0, $$ $$(Z_1,Z_2)\mapsto [Z_1]-[Z_2].$$
By Lemma \ref{lem:3} below, the fiber of this map is a countable union of constructible sets. Thus we can define a dimension $c_d$ \R{of} the general fiber of $\sigma_d$ and set $\dim Im(\sigma_d)=2d\dim U-c_d$.

\begin{definition}
	We say that $h_0(U)^0$ is \emph{infinite-dimensional} if $$\lim_{d\to \infty} \dim \im(\sigma_d)=\infty.$$
\end{definition}

In this paper, using logarithmic algebraic geometry, we prove the log Mumford theorem. 

\begin{thm}[Log Mumford theorem]\label{thm:mumford}
Let $(X,D)$ be a log smooth proper surface pair, and let $U$ be its interior. If $h^0(X,\Omega^2_X(\log D))>0$, then the group $h_0(U)^0$ is infinite-dimensional.
\end{thm}

This is proved in Corollary \ref{cor}. Our proof follows the strategy as in \cite{Mumford68} and the crucial part of the proof is the existence of induced log forms in Proposition \ref{prop:induced}. 

Since the set of $\A^1$-equivalence classes of divisors on open curves is the generalized Jacobian \cite{Serre-AG}, we may formulate the analogue of Bloch's conjecture in the logarithmic setting.

\begin{conj}[Log Bloch conjecture]
Let $(X,D)$ be a log smooth proper surface pair\R{,} and let $U$ be its interior. If $h^0(X,\Omega^2_X(\log D))=0$, then the Albanese morphism induces an isomorphism $$h_0(U)^0\cong \Alb(U).$$ 
\end{conj}

We prove a special case of log Bloch's conjecture as below.

\begin{thm}\label{thm:bloch}
The log Bloch's conjecture holds for log smooth surface pairs with log Kodaira dimension $-\infty$.
\end{thm}

In arbitrary dimension, if $(X,D)$ is log rationally connected, introduced in \cite{CZ,A1,Z3}, then we have the vanishing $h^0(X,\Omega^{\otimes m}_X(\log D))=0$ for any $m$. In this case, we prove that $h_0(U)^0$ vanishes as well.  See Proposition \ref{prop:logRC}. However, this is too weak to prove Theorem \ref{thm:bloch}. There exists an $\A^1$-ruled surface pair with $q(X,D)=\Alb(U)=0$ but not log rationally connected \cite[Section 4]{Z3}.

\begin{notation}
In this paper, we work with (log) varieties and log pairs over complex numbers $\C$. \R{We refer to \cite{KKato} or \cite[Ch. 3]{Gross-b} for basic notions in log geometry.} For any log scheme $(X,\cM_X)$, we denote by $X^\circ$ the open subset with the trivial log structure and denote by $\Omega^q(X,\cM_X)$ the sheaf of log $q$-forms. \R{A \emph{log rational curve} on a log variety $(X,\cM_X)$ is a log morphism $$f:(\P^1,\cM_{\{\infty\}})\to (X,\cM_X),$$
where $\cM_{\{\infty\}}$ is the divisorial log stricture associated to $\{\infty\}$ on $\P^1$. }

A \emph{log pair} $(X,D)$ means a variety $X$ with a reduced Weil divisor $D$. Let $U$ be its interior $X-D$. We say that $(X,D)$ is \emph{log smooth} if $X$ is smooth and $D$ is a normal crossing divisor. A log pair is \emph{proper} if the ambient variety is proper. For a log smooth pair $(X,D)$, we use $\kappa(X,D)$ to denote the \emph{logarithmic Kodaira dimension} and define the \emph{log irregularity} $q(X,D):=h^0(X,\Omega^1_X(\log D))$. Since they only depend on the interior $U$, we may write $\kappa(U)$ and $q(U)$ as well.

\end{notation}

\subsection*{Acknowledgment} 
The author would like to thank Qile Chen, Jason Starr and Burt Totaro for helpful conversations.


\section{Induced log differentials}
Throughout this section, we let $G$ be the symmetric group $S_n$ and let $(X,\cM_X)$ be a log smooth variety over $\C$. Let $D$ be the boundary divisor $X-X^\circ$. By log smoothness, $\cM_X$ is a divisorial log structure 
$$\cM_X=\{f\in\O_X|f\in \O^*_{X-D} \}\subset\O_X.$$

Let $(X^n,\cM_{X^n})$ be the product log structure. Then $\cM_{X^n}$ is $G$-invariant. 

Consider the quotient map:
$$\pi:X^n\to Y:=X^n/G.$$ 
\begin{lemma}\label{lem:G-chart}
	Let $\cM_Y$ be the $G$-invariant subsheaf $\cM_{X^n}^G$. Then $(Y,\cM_Y)$ is a log variety and $$\pi:(X^n,\cM_{X^n})\to (Y,\cM_Y)$$ is a log morphism.
\end{lemma}

\proof Since $(X^n,\cM_{X^n})$ is a log scheme, we have 
$$\O_{X^n}^*\subset \cM_{X^n}\subset \O_X.$$
By taking the $G$-invariant part, we get $$(\O_{X^n}^*)^G\subset \cM_{Y}\subset \O_Y.$$
Since the first term is indeed $\O_Y^*$, we conclude that $Y$ is a log scheme.

The natural diagram
$$\xymatrix{
	\cM_Y \ar[d]\ar[r] &\O_{Y}\ar[d]^{\pi^*}\\
	\cM_{X^n} \ar[r] &\O_X},$$
where all arrows are inclusions shows that $\pi$ is a log morphism.
\qed

\begin{lemma}\label{lem:Yfs}
	\R{The log variety} $(Y,\cM_Y)$ is fine and saturated. 
\end{lemma}

\proof We know that \'etale locally on $X$, there exists a fine and saturated chart $$P\to \O_X.$$ Furthermore, by choosing the defining equations of the irreducible components of $D$, we may assume the chart morphism factors as below:
$$P\subset \cM_X\subset \O_X,$$
and $\cM_X$ is isomorphic to $P\oplus \O_X^*$.
This induces a $G$-invariant fs chart 
$$P^n\to \O_{X^n}$$ for $(X^n,\cM_{X^n})$ such that $$\cM_{X^n}\cong P^n\oplus \O_{X^n}^*.$$
Now taking the $G$-invariant part, we get $$\cM_Y\cong P\oplus \O_Y^*,$$ and actually $P$ maps to the defining equations of the boundary divisors on $Y$. Therefore, $(Y,\cM_Y)$ is a fine saturated log scheme.\qed

\ \\

For any log smooth variety $(S,\cM_S)$ with a morphism $f:(S,\cM_S)\to (Y,\cM_Y)$, let \R{${S'}=(S\times_Y X^n)$} be the fibered product with the fs log structure $\cM_{{S}'}$, c.f. \cite[II.2.4]{Ogus}. Let \R{$\tilde{S}=(S\times_Y X^n)_{red}$} with the induced log structure from $S'$. We have a diagram as below:
$$\xymatrix{
(\tilde{S},\cM_{\tilde{S}}) \ar[rd]^p \ar[r]^i &({S'},\cM_{{S'}}) \ar[d]^{p'}\ar[r]^{\tilde{f}} &(X^n,\cM_{X^n})\ar[d]^\pi\\
&(S,\cM_S) \ar[r]^f &(Y,\cM_Y)}$$

Given a $G$-invariant log $q$-form $\omega\in \Gamma(X^n,\Omega^q(X^n,\cM_{X^n}))$, let 
$$\tilde\omega=(\tilde{f}\circ i)^*(\omega)\in \Gamma(\tilde{S}, \Omega^q(\tilde{S},\cM_{\tilde{S}})).$$
Then $\tilde{\omega}$ is $G$-invariant.

\begin{proposition}\label{prop:induced}
If $S$ is log smooth, there exists a unique log $q$-form $\eta_f\in \Gamma({S}, \Omega^q({S},\cM_{{S}}))$ such that 
$$p^*(\eta_f)-\tilde{\omega} \text{ is torsion in } \Omega^q(\tilde{S},\cM_{\tilde{S}}).$$
\end{proposition}

\begin{remark}\label{rem}
	When $S$ has the trivial log structure, this construction of $\eta_f$ coincides with the construction in \cite[Section 1]{Mumford68}.
\end{remark}

\proof First we prove the uniqueness. Indeed, there are non-singular open dense subsets $S_0\subset S$, $\tilde{S_0}=p^{-1}(S_0)\subset \tilde{S}$ with trivial log structures such that \R{$S_0=\tilde{S_0}/\Gamma$ and $\Gamma$ acts freely on $\tilde{S_0}$, where $\Gamma$ is a quotient group of the stabilizer group of the open subset $\tilde{S_0}$}. Thus $\tilde\omega|_{\tilde{S_0}}$ as a regular form \R{descends} to a regular form $\theta$ on $S_0$. By the condition in the Lemma, $\eta_f$ coincides with $\theta$ over $S_0$, thus is unique.

Let $\eta_f$ be the meromorphic form extending $\theta$ on $S$. To prove the existence, it suffices to check that $\eta_f$ as a meromorphic section of $\Omega^q(S,\cM_S)$ is regular everywhere. Since $(S,\cM_S)$ is log smooth, hence $S$ is normal, it suffices to check this at points of codimension one. Hence we may assume that $S$ is the spectrum of a local discrete valuation ring $R$ with the fraction field $K$. Let $T$ be the normalization of $\tilde{S}$ and consider the normalization morphism $$a:T\to \tilde{S}.$$ The morphism $p'$ is finite, so is the composite morphism $$p\circ a: T\to S.$$ In particular, $T$ is a disjoint union of local discrete valuation ring $T_i=\Spec R_i$ with the generic point $\Spec K_i$. The log structure on $T$ is given canonically below. 
\begin{lemma}\label{lem:canonical}
	There exists a canonical fs log structure on $T$ by choosing $$\cM_{T_i}=R_i-{0}\subset \O_{T_i}=R_i.$$ In particular, $(T,\cM_T)$ is log smooth. \qed
\end{lemma}

\begin{lemma}
	The morphism $a:T\to \tilde{S}$ extends to a unique log morphism:
	$$a:(T,\cM_T)\to (\tilde{S},\cM_{\tilde{S}}).$$
\end{lemma}

\proof We may assume that both $T$ and $\tilde{S}$ are irreducible. Since $({S'},\cM_{{S'}})$ is fine and saturated, there exists a fs chart 
$$c:P\to \O_{{S'}}.$$ 
To show that $a$ is a log morphism, it suffices to prove the image of the composite morphism 
$$P\to \O_{S'}\to i_*\O_{\tilde{S}}\to (i\circ a)_*\O_T$$
does not contain zero. Since $a$ is the normalization map, it is enough to show the image of $P$ in $\O_{\tilde{S}}$ does not contain zero, or equivalently, none of the \R{images} of $P$ in $\O_{S'}$ is nilpotent.

If there exists $p\in P$ such that $c(P)$ is nilpotent, then consider the base change of $c(p)\otimes_{\O_{S}}K$ via the following diagram is still nilpotent.

$$\xymatrix{
	\O_{S'}\otimes_R K &\O_{S'}\ar[l] \\
	K \ar[u] &\O_S=R \ar[l]\ar[u]^{p'^*}}$$
Indeed, we have that $\O_{S'}$ is a flat $\O_S$-module and $\O_S$ is a principal ideal domain. Thus $\O_{S'}$ is torsion free. In particular, the nilpotent elements cannot be killed after tensoring with $K$.

This tells us the log structure on $S'$ is nontrivial over $\Spec K$. On the other hand, since $S^\circ$ is nonempty, we have a log morphism 
$$(\Spec K, \text{ trivial log structure}) \to (S,\cM_S)$$
which induces a Cartesian diagram
$$\xymatrix{
	S'\otimes_S \Spec K \ar[d]\ar[r] &(X^n,\cM_{X^n})\ar[d]^{}\\
	(\Spec K, \text{ trivial log structure}) \ar[r] &(Y,\cM_Y)}.$$

By the universal property of log fibered product, $S'\otimes_S \Spec K$ must have the trivial log structure. \R{This is a contradiction}.\qed

 Now let us return to the proof of Proposition \ref{prop:induced}, we construct a diagram
$$\xymatrix{
	(T,\cM_T)\ar[r]^a \ar[rrd]_{r=p\circ a} &(\tilde{S},\cM_{\tilde{S}}) \ar[rd]^p \ar[r]^i &({S'},\cM_{{S'}}) \ar[d]^{p'}\ar[r]^{\tilde{f}} &(X^n,\cM_{X^n})\ar[d]^\pi\\
&	&(S,\cM_S) \ar[r]^f &(Y,\cM_Y)}$$
such that 
\begin{itemize}
	\item $(T,\cM_T)$ is log smooth;
	\item $r$ is finite.
\end{itemize}

Since $p^*(\eta_f)-\tilde{\omega}$ is torsion and $(T,\cM_T)$ is log smooth, we have 
$$r^*(\eta_f)=a^*(p^*(\eta_f))=a^*(\tilde{\omega}).$$
Since $\omega$ as an element in $\Gamma(X^n,\Omega^1(X^n,\cM_{X^n}))$ is regular, $$r^*(\eta_f)=a^*(\tilde{\omega})=(\tilde{f}\circ i\circ a)^*\omega$$ is a regular as an element in $\Gamma(T,\Omega^q(T,\cM_{T}))$. Now the proposition is proved using the following lemma.\qed

\begin{lemma}
	There is a well defined trace map$$tr: \Omega^1(T,\cM_T)\to \Omega^1(S,\cM_S)$$ such that the composite	
	$$\xymatrix{
		\Omega^1(S,\cM_S) \ar[r]^-{r^*} & \Omega^1(T,\cM_T) \ar[r]^{tr} &\Omega^1(S,\cM_S) }$$
	is multiplication by the degree of $r$.
\end{lemma}

\proof By construction, if $(S,\cM_S)$ has the trivial log structure, so does $(T,\cM_T)$. Now we can simply use the standard trace map, c.f., Mumford's paper \cite{Mumford68}. From now on, we assume $(S,\cM_S)$ has nontrivial log structure, so does $(T,\cM_T)$. Furthermore, since $S$ is the spec of a local ring and log smooth, the log structure $\cM_S$ is the canonical one as in Lemma \ref{lem:canonical}. Let $\mathfrak{m}_S$, $\mathfrak{m}_T$ be the maximal ideals respectively.

We claim that the morphism $$r:(T,\cM_T)\to (S,\cM_S)$$ is log \'etale. Consider the commutative diagram given by the charts
$$\xymatrix{
	(T,\cM_T)  \ar[d]\ar[r] &\A^1=\Spec \Z[\N]\ar[d]^{u}\\
	(S,\cM_S) \ar[r] &\A^1=\Spec \Z[\N]}.$$
Here the map $u$ is $t\mapsto t^k$, where $\mathfrak{m}_S \O_T=\mathfrak{m}_T^k.$ This implies that the natural morphism 
$$T\to S\times_{\A^1} \A^1$$ is unramified. Therefore, $r$ is log \'etale.

By the universal properties of log differentials, we have a sequence
$$\Omega^1(S,\cM_S)\otimes_{\O_S}\O_T \to \Omega^1(T,\cM_T)\to \Omega^1_{(T,\cM_T)|(S,\cM_S)}\to 0.$$
Since $r$ is log \'etale, the last term vanishes. Since both $(T,\cM_T)$, $(S,\cM_S)$ are log smooth of dimension one, by Nakayama Lemma, we have the isomorphism $$r^*:\Omega^1(S,\cM_S)\otimes_{\O_S}\O_T \to \Omega^1(T,\cM_T)$$

Now the log trace map is constructed as below:
	$$\xymatrix{
		\Omega^1(T,\cM_T) \ar[r]^>>>>{(r^*)^{-1}} & \Omega^1(S,\cM_S)\otimes_{\O_S}\O_T \ar[r]^-{tr} &\Omega^1(S,\cM_S) },$$
	where the second map is induced by the trace map $tr:\O_T\to \O_S$. The second part of the Lemma trivially follows.\qed

\section{Log Mumford's theorem}

\begin{lemma}\label{lem:extendslog}
\R{	Given a proper log variety $(V,\cM_V)$ and a normal scheme $T$, any morphism $$\A^1\times T\to V^\circ$$ uniquely extends to a family of log rational curves over $T_0$ 
	$$(\P^1,\cM_{\{\infty\}})\times T_0\to (V,\cM_V),$$ where $T_0$ is a dense open subset of $T$ and $\cM_{\{\infty\}}$ is the divisorial log structure associated to ${\{\infty\}}\subset \P^1$.}
\end{lemma}

\proof Since $V$ is proper and $T$ is normal, we have a morphism $$u:\P^1\times T_0\to V,$$ \R{where $T_0$ is a dense open subset of $T$. Consider the commutative diagram
	$$\xymatrix{
	u^{-1}\cM_V  \ar[d]^\alpha & \cM_{\{\infty \}\times T_0} \ar@{^{(}->}[d]^i \\
	u^{-1}\O_V \ar[r]^{u^*} &\O_{\P^1\times T_0}
 }.$$
To prove $u$ extends to a log morphism, it suffices to prove that for any element $g\in \cM_D$, $u^*(\alpha(g))$ lies in $\cM_{\{\infty \}\times T_0}\subset u_*\O_{\P^1\times T_0}$, or equivalently, $u^*(\alpha(g))$ is invertible on $\A^1\times T_0$. By assumption, the image of $\A^1\times T_0$ under $u$ factors through $V^\circ$. Thus we have $$u^*(\alpha(g))|_{\A^1\times T_0}=u^*(\alpha(g)|_{V^\circ})|_{\A^1\times T_0}.$$
Since the log structure on $V^\circ$ is $\O_{V^\circ}^*$, we have $\alpha(g)|_{V^\circ}\in \O_{V^\circ}^*$. In particular, $u^*(\alpha(g))$ is invertible on $\A^1\times T_0$.
}\qed

\begin{notation}
	Let $(X,D)$ be a log smooth \R{proper variety} with the interior $U$. Let $G=S_n$. We pick a nonzero logarithmic \R{$q$-form $\omega\in \Gamma(X,\Omega^q_X(\log D))$}. Let $\omega^{(n)}=\sum_{1}^n p_i^*\omega\in \Gamma(X^n,\Omega^q(X^n,\cM_{X^n}))$. Then $\omega^{(n)}$ is $G$-invariant. By Proposition \ref{prop:induced}, for \R{every} log smooth variety $(S,\cM_S)$ and morphism 
	$$f:(S,\cM_S)\to (Y,\cM_Y),$$ we have an \R{induced $q$-form $$\eta_f\in \Gamma(S,\Omega^q(S,\cM_S)).$$}
\end{notation}

\begin{thm}\label{thm:eta0}
	Let $T$ be a smooth variety. Given a morphism $f:T\to S^nU$, it extends to a morphism 
	$$f:(T,\O_T^*)\to (Y^n,\cM_{Y^n}).$$ If all the $0$-cycles in the image $f(T)$ are $\A^1$-equivalent, then it follows that $$\eta_f=0.$$
\end{thm}

\begin{lemma}\label{lem:3}
	$S^nU\times S^n U$ contains a countable set $Z_1,Z_2,\cdots$ of constructible sets, such that if $(A,B)\in S^nU\times S^nU$, then 
	$$A\sim_{\A^1} B\iff (A,B)\in \cup_{i=1}^\infty Z_i.$$
	For each $i$, there is a reduced scheme $W_i$ and a set of morphisms
	$$e_i:W_i\to Z_i,$$
	$$f_i:W_i\to S^m U,$$
	$$g_i:W\times \A^1\to S^{n+m}U$$ such that we get the equations between zero-cycles:
	$$g_i(w,0)=p_1(e_i(w))+f_i(w),$$
	$$g_i(w,1)= p_2(e_i(w))+f_i(w),$$
 for all $w\in W_i$ and $e_i$ is surjective.
\end{lemma}

\proof We observe the fact that if $A,B\in S^k U$ are joined by a chain of $p$ $\A^1$-curves $E_1,\cdots,E_p$ such that $$E_1(0)=A,$$ $$E_{p}(1)=B,$$
$$E_i(1)=E_{i+1}(0)=C_i,i=1,\cdots,p,$$then $A+C_1+\cdots+C_{p-1}$ and $C_1+\cdots+C_{p}+B$ in $S^{pk}U$ are joined by a single $\A^1$-curve, whose degree is bounded by the degree of the $E_i$'s. Therefore, for any pair $(A,B)$, the condition $A\sim_{\A^1} B$ is equivalent to that there exists $C\in S^mU$ and an irreducible $\A^1$-curve $E$ on $S^{n+m}U$ of bounded degree connecting $A+C$ and $B+C$.

For any $l$, we define $(Y^l,\cM_{Y^l})$ the fine saturated log scheme constructed in Lemma \ref{lem:G-chart} and Lemma \ref{lem:Yfs} for the quotient scheme $Y^l:=X^l/S_l$. Clearly, there exists \R{a} strict open immersion 
$$(S^lU,\O_{S^lU}^*)\to(Y^l,\cM_{Y^l}). $$ By Lemma \ref{lem:Yfs} and Lemma \ref{lem:extendslog}, any $\A^1$-curve on $S^lU$ extends uniquely to a \R{log rational curve} on $(Y^l,\cM_{Y^l}). $

Now let $\cA_2(Y^{n+m},\cM_{Y^{n+m}};{\le p})$ be the moduli space of two-pointed stable \R{log rational curves} of degree $\le p$ on $(Y^{n+m},\cM_{{n+m}})$, c.f., \cite{GS, Chen,AC}  and let $$\cA^{\circ}_2(Y^{n+m},\cM_{Y^{n+m}};{\le p})\subset \cA_2(Y^{n+m},\cM_{Y^{n+m}};{\le p})$$ be the log trivial part which parametrize two-pointed \R{log rational} curves. We have the natural evaluation morphism $$ev_{{n+m},p}:\cA^{\circ}_2(Y^{n+m},\cM_{Y^{n+m}};{\le p})\to S^{m+n}U\times S^{m+n} U$$ Define the incidence reduced subscheme
$$W_{{n+m},p}\subset S^nU\times S^nU\times S^mU\times \cA^{\circ}_2(Y^{n+m},\cM_{Y^{n+m}};{\le p}),$$
$$W_{{n+m},p}=\{((A,B), C,g)|g(0)=A+C,g(1)=B+C \}.$$
Define $Z_{{n+m},p}$ as the image of $W_{{n+m},p}$ under the projection to $S^nU\times S^nU$, which is constructible. Define $e_{n+m,p}$, $f_{n+m,p}$ the restriction of the natural projection morphisms on $W_{n+m,p}$. The morphism $g_{n+m,p}$ are defined via the universal morphism of log rational curves on $\cA^{\circ}_2(Y^{n+m},\cM_{Y^{n+m}};{\le p})$.  \qed

\begin{remark}
	\R{In the proof of Lemma \ref{lem:3}, the moduli space of log rational curves are not really needed. Any other reasonable sequence of moduli spaces could also be used to define the constructible sets $\{Z_i\}$, for example, \cite[Def. 5.1, Prop. 5.3]{KM}.
}
\end{remark}

\proof[Proof of Theorem \ref{thm:eta0}] Given $f:S\to S^nU$ such that all zero cycles $f(s)$ are $\A^1$-equivalent, fix a base point $A_0$ in the image. It follows from Lemma \ref{lem:3} and Lemma \ref{lem:extendslog} that there is a non-singular variety $T$, a dominant morphism $e:T\to S$ and morphisms $$g: T\to S^mU,$$ $$h:(\P^1,\cM_{\{\infty\}})\times T\to (Y^{n+m},\cM_{Y^{n+m}}),$$ such that:
$$h(t,0)=g(t)+f(e(t)),$$ $$h(t,1)=g(t)+A_0,$$ for all $t\in T$. 

By Proposition \ref{prop:induced} and Lemma \ref{lem:extendslog}, we have induced log \R{$q$}-forms $\eta_f$, $\eta_g$ and $\eta_h$. By Remark \ref{rem}, we note that $\eta_f$, $\eta_g$, $\eta_h|_{T\times \{0\}}$ and $\eta_h|_{T\times \{\infty\}}$ are indeed regular \R{$q$-forms} constructed by Mumford. By \cite[Lemma 2]{Mumford68}, we have 
$$\eta_h|_{T\times \{0\}}=\eta_g+e^*(\eta_f),$$
$$\eta_h|_{T\times \{\infty\}}=\eta_g+\eta_{A_0}.$$

Now $\eta_h$ is a \R{log $q$-form} on $(\P^1,\cM_{\{\infty\}})\times T$. Since
$$\Omega^q{((\P^1,\cM_{\{\infty\}})\times T)}\cong p_1^*(\Omega^q_T)+p_1^*(\Omega_T^{q-1})\otimes p^*_2(\Omega^1(\P^1, \cM_{\{\infty \}}))$$
and $\Omega^1(\P^1, \cM_{\{\infty \}})\cong \O_{\P^1}(-1)$ has no global sections, it follows that $$\eta_h=p_1^*(\eta)$$ for some $\eta\in \Gamma(\Omega^q_T)$. Therefore, 
$$\eta_h|_{T\times \{0\}}=\eta_h|_{T\times \{\infty \}}.$$ Since $\eta_{A_0}=0$, we find $e^*(\eta_f)=0$, hence $\eta_f=0.$\qed 

\ \\
\R{Now let us assume that $\dim U$=2 and $q=2$. }Let $(S^nU)_0$ be the open subset parametrizing zero-cycles $\sum_{i=1}^n x_i$ such that $x_i$'s are all distinct and $\omega(x_i)\neq 0$ for all $i$. The open immersion $$f: (S^nU)_0\to S^nU$$ induces a log morphism $$f:(S^nU)_0\to (Y^n,\cM_{Y^l}).$$ The induced log $2$-form is \R{a} holomorphic symplectic form. The maximal isotropic subspace of $\eta_f$ is of dimension $n$. If $S\subset (S^nU)_0$ is a nonsingular subvariety parametrizing $\A^1$-equivalent zero-cycles, we have $\eta_f|_S=0$, thus $\dim S\le n$.

\begin{corollary}\label{cor}
	Let $(X,D)$ be a log smooth surface with $h^0(\Omega^2_X(\log D))> 0$ and let $(S^nU)_0$ be defined as above. Then if $S\subset (S^nU)_0$ is a subvariety consisting $\A^1$-equivalent $0$-cycles, it follows that $\dim S\le n$.\qed
\end{corollary}

\section{The Log Bloch conjecture}

\subsection{Log rationally connected varieties}

\begin{lemma}\label{lem:move}
	Let $U$ be a smooth quasiprojective curve. For any \R{dense} open subset $V\subset U$, any point $x\in U$ is $\A^1$-equivalent to $A-B$, where both $A$ and $B$ are effective divisors supported on $V$.
\end{lemma}

\proof If $x\in V$, then the lemma is trivial. We assume that $x\notin V$. We choose the compactification $(X,D)$ of $U$ with $D=p_1+\cdots+p_d$, where all $p_i$'s are distinct.  We pick an effective divisor $B\subset V$ with sufficiently high degree satisfying
\begin{itemize}
	\item $h^0(\O(x+B-D))=h^0(\O(x+B))-d$;
	\item $\O(x+B-D)$ is very ample.
\end{itemize}

Let $H_i$ be the hyperplane in $|x+B|$ parametrizing divisors containing $p_i$. By the above condition, \R{the} $H_i$'s intersect transversally in $|x+B|$ and a divisor in $|x+B|$ is away from $D$ if and only if it avoids $\cup_{i=1}^d H_i$. Since $\O(x+B-D)$ is very ample, we may choose an effective divisor $E\in H^0(\O(x+B-D))$ and $E\subset V\backslash B$. The base-point-free pencil connecting $x+B$ and $D+E$ is an $\A^1$-curve on the pair $(|x+B|, \cup_{i=1}^d H_i)$. Let $A$ be a general element of this pencil. Then $A\sim_{\A^1} x+B$. Since $B\subset V$ and the pencil is base point free, $A$ is supported on $V$ as well.\qed

\begin{lemma}\label{lem:open}
	Let $U$ be a smooth quasiprojective variety and $V\subset U$ be a \R{dense} open subset. Then the natural map
	$$i_*:h_0(V)\to h_0(U)$$ is surjective.
\end{lemma}

\proof By choosing a smooth curve $C$ on $U$ with $x\in C$ and $C\cap V\neq\emptyset$, Lemma \ref{lem:move} implies that any point $x\in U-V$ is $\A^1$-equivalent to $A-B$, where both $A$ and $B$ are effective zero cycles on $V$. The lemma follows.\qed

\begin{proposition}\label{prop:logRC}
	If $(X,D)$ is log rationally connected, then \R{$h_0(U)=\Z$}.
\end{proposition}

\proof Since $(X,D)$ is log rationally connected, let $p$ be a general point on $U$ and let $ U'\subset U$ be a nonempty open subset of $U$ such that any point in $U'$ is connected by a log rational curve through $p$. Thus $h_0(U')=\Z$. \R{By Lemma \ref{lem:open}, $h_0(U)$ is isomorphic to $\Z$ as well}.\qed  

\begin{remark}
	In general, we do not know that any pair of points in the interior of a log RC pair is connected by a log rational curve. Any pairs with such properties are called \emph{strongly log RC pairs}.
\end{remark}

\subsection{Surface pairs with $\kappa=-\infty$}

\proof[Proof of Theorem \ref{thm:bloch}] 
Let $(X,D)$ be a proper log smooth surface pair with $\kappa(X,D)=-\infty$. By \cite[Theorem 3.47]{Kollar-Mori}, we run the log minimal model program on this pair
$$(X,D)=(X_0,D_0)\to(X_1,D_1)\to\cdots\to (X_k,D_k)=(X^*,D^*),$$
such that: 
\begin{enumerate}
	\item the log Kodaira dimension remains the same, i.e., $\kappa(X_i,D_i)=-\infty$;
	\item the end product $(X^*,D^*)$ is either \begin{enumerate}
		\item log ruled, or 
		\item a log del Pezzo surface of Picard number one, i.e., $\rho(X^*)=1$.
	\end{enumerate}
\end{enumerate}

If the minimal model $(X^*,D^*)$ is log del-Pezzo surface but not log ruled, then by the works of Miyanishi-Tsunoda \cite{MT84a}\R{,} Keel-McKernan \cite{KM} and \cite[Lemma 2.1, Theorem 2.2, 2.3]{Z3}, $(X,D)$ is log rationally connected. In this case, Theorem \ref{thm:bloch} follows from Proposition \ref{prop:logRC}. 

If the minimal model $(X^*,D^*)$ is log ruled, then by \cite[Lemma 2.1]{Z3}, $(X,D)$ is log ruled. In this case, Theorem \ref{thm:bloch} follows from Proposition \ref{prop:0} below.\qed

\begin{proposition}\label{prop:0}
Log Bloch's conjecture holds for log ruled surface pairs.
\end{proposition}

First we observe the following lemma.
\begin{lemma}\label{lem:blow}
Let $(X,D)$ be a log smooth proper surface pair with the interior $U$. Let $X'$ be the surface obtained by a sequence of blow ups on $X$:
$$b:X'\to X,$$ with the boundary $D':=b^{-1}(D)$. Then log Bloch's conjecture holds for $(X,D)$ if and only if it holds for $(X',D')$.
\end{lemma}
\proof Let $U'$ be the interior of $(X',D')$. We have a commutative diagram as below
$$\xymatrix{
	h_0(U')^0\ar[d]_{b_*}\ar[r] &\Alb(U')\ar[d]^{\Alb(b)}\\
	h_0(U)^0 \ar[r] &\Alb(U).}$$
Since blowing up does not change the Albanese, it suffices to show that $$b_*:h_0(U')^0\to h_0(U)^0$$ is an isomorphism. This follows from the blowing up long exact sequence of Suslin's algebraic singular homology \cite[Proposition 14.19]{Veovo-MC}.\qed

\proof[Proof of Proposition \ref{prop:0}] Lemma \ref{lem:blow} implies that\R{,} without loss of generality, we may always replace $(X,D)$ by a sequence of blow ups to prove Proposition \ref{prop:0}. Now assume that $(X,D)$ is log ruled surface. Let $q$ be the log irregularity $h^0(\Omega^1_X(\log D))$. We first construct a log ruling on $U$ based on the value of $q$. {\ \\}
	
	{\bf The case when $q>0$}  
	
	Consider the Albanese morphism \cite{Iitaka76}
	$$a:U\to \Alb(U),$$
	where $\Alb(U)=H^0(\Omega^1_X(\log D))^*/H_1(U,\Z)$ as a semiabelian variety of dimension $q$. Let $T_0$ be the closure of the image $a(U)$ and we rename the map $a:U\to T_0$ as $$f:U\to T_0.$$
	Since $h^0(K_X+D)=0$, $T_0$ is a curve on $\Alb(U)$. Otherwise, any nowhere vanishing log $2$-form on $\Alb(U)$ pulls back to a nonzero log $2$-form on $U$. By \cite[Corollary 1]{Iitaka76}, $T_0$ is a smooth curve with the diagram below 
	$$\xymatrix{
		U\ar[d]_f\ar[r]^-a &\Alb(U)\ar[d]_{\Alb(f)}^{\cong}\\
		T_0 \ar[r] &\Alb(T_0),}$$ and a general fiber of $f$ is irreducible.
	
	\begin{lemma}
		The morphism $f:U\to T_0$ is surjective and it gives the log ruling on $U$, that is, a general fiber of $f$ is a log rational curve.
	\end{lemma}
	
	\proof Since there are no log rational curves on the Albanese, the log ruling on $U$ gets contracted via $f$. Since the general fiber of $f$ is irreducible, $f$ gives the log ruling. Denote the image $f(U)$ by $T'_0\subset T_0$. Since every log one form on $T'_0$ pulls back to a log one form on $U$ and $q(T_0)\le q(T'_0)$, we have 
	$$q(T_0)\le q(T'_0)\le q(U)=q(T_0).$$
	This implies $q(T'_0)=q(T_0)$. Thus $T'_0=T_0$.\qed	
	 {\ \\}
	 
	 {\bf The case when $q=0$}  
	 
	 After blowing up finitely many points on $U$, still denoted by $U$, we pick a log ruling: $$f:U\to T_0,$$ where $T_0$ is a smooth curve. We may further assume $f$ is surjective. 
	 \begin{lemma}
	 	$T_0$ is either $\P^1$ or $\A^1$.
	 \end{lemma}
	 \proof Pullback of log one forms under $f$ implies $$q(T_0)\le q(U)=0. $$\qed
	
In either case, we choose the smooth compactification $T$ of $T_0$ and let $S=T-T_0$. After further blowing up finitely many points, still \R{denoting} the pair $(X,D)$ and the interior $U$, we may assume there exists a proper flat morphism $$f:(X,D)\to (T,S),$$ such that\begin{enumerate}
	\item $f(U)=T_0$; and 
	\item the morphism $$\Alb(f):\Alb(U)\to \Alb(T_0)$$ is an isomorphism. 
\end{enumerate} 

Let $T_1$ be an open subset of $T_0$ and let $U_1:=f|_U^{-1}(T_1)$ be an open subset of $U$ such that $$U_1\cong T_1\times \A^1.$$Choose a regular section $\sigma:T_1\to U_1$. We have the diagram as below:
 $$\xymatrix{
	T_1\ar[rd]_{j}\ar[r]^\sigma &U\ar[d]^{f}\\
	 &T_0,}$$
where \R{$j$} is indeed an open immersion. It induces a diagram on algebraic singular homology $$\xymatrix{
	h_0(T_1)^0\ar[rd]_{j_*}\ar[r]^{\sigma_*} &h_0(U)^0\ar[d]^{f_*}\\
	&h_0(T_0)^0.}$$
\begin{lemma}\label{lem:curve}	
		\R{We have a commutative diagram as below
		$$\xymatrix{
			h_0(T_1)^0\ar[d]_{j_*}\ar[r]^{a_1} &\Alb(U)\ar[d]_{\Alb(j)}\\
		h_0(T_0)^0\ar[r]^{a_0}	&\Alb(T_0)}$$
	where both $a_0$ and $a_1$ are isomorphisms. Furthermore, $j_*$ is surjective with the kernel an algebraic torus $H$.}

\end{lemma}
\proof It follows from the theory of generalized Jacobians \cite{Serre-AG}.\qed 

\begin{lemma}\label{lem:j}
	$\sigma_*$ is surjective. 
\end{lemma}
\proof The \R{inclusion map $i:\sigma(T_1)\to U$} factors as below:
$$\xymatrix{
	\sigma(T_1)\ar[r]^-{i_1} & U_1=T_1\times \A^1\ar[r]^-{i_2}& U}.$$
The lemma follows from that  $i_{1*}$ induces an isomorphism on $h_0$ and that $i_{2*}$ is surjective by Lemma \ref{lem:open}.\qed

By Lemma \ref{lem:curve} and \ref{lem:j}, we have the following commutative diagram:
 $$\xymatrix{
 \Alb(T_1)\cong	h_0(T_1)^0\ar@{->>}[rd]_-{j_*}\ar@{->>}[r]^-{\sigma_*} & h_0(U)^0\ar[d]^{f_*}\ar[r]& \Alb(U)\ar[d]_\cong^{\Alb(f)}\\
 	&h_0(T_0)^0\ar[r]_\cong & \Alb(T_0).}$$
 Proposition \ref{prop:0} follows from the following lemma.\qed

\begin{lemma}
$$\ker(\sigma_*)=\ker(j_*).$$
\end{lemma}

\proof We only need to prove $\ker(j_*)\subset ker(\sigma_*)$. By Lemma \ref{lem:curve}, any element in $\ker(j_*)$ is of the form $A-B$ satisfying 
\begin{itemize}
	\item both $A$ and $B$ are effective divisors on $T_1$ of degree $d$;
	\item $A\sim_{\A^1} B$ on $T_0$.
\end{itemize}
Note that if $T_0=T$ is proper, we may replace $S=\emptyset $ by $S=\{p \}$, where $p$ is away from the support of $A+B$. Then $A$ is $\A^1$-equivalent to $ B$ on the open curve $T_0-\{p\}$ as well. So for the rest of the proof, we assume that $S\neq\emptyset$.

Recall that the morphism $$f: (X,D)\to (T,S)$$ is log ruled. By reduced fiber theorem \cite[\href{http://stacks.math.columbia.edu/tag/09IL}{Tag 09IL}]{stack}, after a finite base change \R{and also normalizing}, 
$$
\xymatrix{
	(X',D') \ar[r]^{g'} \ar[d]^{f'} & (X,D)\ar[d]^f \\
	(T',S'=g^{-1}(S)) \ar[r]^-g & (T,S),
}
$$
we may assume that all geometric fibers of \R{$f'$} are reduced over $T'-S'$. Since $(X',D')$ is log ruled over $(T',S')$, strong approximation over complex function fields away from $S'\neq\emptyset$ \cite[Theorem 6.13]{Rosen2002} implies that there exists an integral section $s'$ over $T'-S'$. The image $R_0:=g'(s'(T'-S'))$ gives a integral multisection of $f$ which is finite of degree $N$ over $T_0$. Let $u:R_0\to T_0$ be the natural map. 

Since Suslin's homology is contravariant for finite flat maps \cite[Section 4]{Geisser2010}, the equivalence $A\sim_{\A^1} B$ on $T_0$ implies that $$u^*(A)\sim_{\A^1} u^*(B)$$ on $R_0$. On the other hand, by our construction, for any $p\in T_1$, $f^{-1}(p)$ is a log rational curve, in particular, $u^*(p)\sim_{\A^1} N\sigma(p)$. Thus we have $$u^*(A)\sim_{\A^1} N\sigma_*(A)$$ $$u^*(B)\sim_{\A^1} N\sigma_*(B).$$ It follows that $$N(\sigma_*(A-B))\sim_{\A^1} 0.$$ Since $\sigma_*(A-B)$ is torsion in $h_0(U)^0$ and it maps to zero under $f_*$, by \cite[Theorem 1.1]{Szamuely03}, $\sigma_*(A-B)$ is trivial on $h_0(U)^0$.\qed

\bibliographystyle{amsalpha}             
\bibliography{myref}

\end{document}